\documentclass[aos,MSNbibl,seceqn,nameyear,dvips]{arximspdf}
\usepackage{dcolumn}
\usepackage{graphicx}

\doi{10.1214/11-AOS902}
\volume{39}
\issue{5}
\pubyear{2011}
\firstpage{2356}
\lastpage{2382}

\makeatletter

\def\bsuffix #1{#1}

\newcolumntype{d}[1]{D{.}{.}{#1}}

\newtheorem{theorem}{Theorem}
\newproclaim{scen}{Scenario}

\newcommand{\F}{\mathcal{F}}
\newcommand{\g}{\mathcal{G}}
\newcommand{\hp}{\hat{p}}
\newcommand{\hl}{\hat{L}}
\newcommand{\till}{\tilde{L}}
\newcommand{\tsig}{\tilde\sigma}
\newcommand{\tils}{\tilde{s}}
\newcommand{\hbp}{\hat{\mathbf{p}}}
\newcommand{\hv}{\hat{v}}
\newcommand{\barp}{\bar{p}}
\newcommand{\bary}{\bar{Y}}
\newcommand{\hmu}{\hat\mu}
\newcommand{\barm}{\bar\mu}
\newcommand{\hf}{\hat{f}}
\newcommand{\hY}{\hat{Y}}
\newcommand{\hdel}{\hat\Delta}
\newcommand{\bx}{\mathbf{x}}
\newcommand{\bbet}{{\bolds\beta}}

\makeatother

\begin{document}
\begin{frontmatter}

\title{Evaluating probability forecasts}
\runtitle{Evaluating probability forecasts}

\begin{aug}
\author[A]{\fnms{Tze Leung} \snm{Lai}\corref{}\thanksref{t1}\ead[label=e1]{lait@stanford.edu}},
\author[B]{\fnms{Shulamith T.} \snm{Gross}\thanksref{t2}\ead[label=e2]{Shulamith.Gross@baruch.cuny.edu}}
\and
\author[C]{\fnms{David Bo} \snm{Shen}\ead[label=e3]{happybo7@gmail.com}}
\runauthor{T. L. Lai, S. T. Gross and D. B. Shen}
\affiliation{Stanford University, Baruch College/CUNY and UBS}
\address[A]{T. L. Lai\\
Department of Statistics\\
Stanford University\\
Sequoia Hall, 390 Serra Mall\\
Stanford, California 94305-4065\\
USA\\
\printead{e1}} 
\address[B]{S. T. Gross\\
Zicklin School of Business\\
Baruch College\\
City University of New York\\
Box 11-220, 1 Bernard M. Baruch Way\\
New York, New York 10010\\
USA\\
\printead{e2}}
\address[C]{D. B. Shen\\
UBS\\
677 Washington Boulevard\\
Stamford, Connecticut 06901\\
USA\\
\printead{e3}}
\end{aug}

\thankstext{t1}{Supported in part by NSF Grant DMS-08-05879.}
\thankstext{t2}{Supported in part by PSC-CUNY 2008 and 2009 grants and
a 2008 Summer Research Support grant from Baruch Zicklin School of
Business.}

\received{\smonth{2} \syear{2011}}
\revised{\smonth{5} \syear{2011}}

%
\begin{abstract}
Probability forecasts of events are routinely used in climate
predictions, in forecasting default probabilities on bank loans or
in estimating the probability of a patient's positive response to
treatment. Scoring rules have long been used to assess the efficacy
of the forecast probabilities after observing the occurrence, or
nonoccurrence, of the predicted events. We develop herein a
statistical theory for scoring rules and propose an alternative
approach to the evaluation of probability forecasts. This approach
uses loss functions relating the predicted to the actual
probabilities of the events and applies martingale theory to
exploit the temporal structure between the forecast and the
subsequent occurrence or nonoccurrence of the event.
\end{abstract}

%

\end{frontmatter}

\section{Introduction}\label{sec1}

Probability forecasts of future events are widely used in diverse
fields of application. Oncologists routinely predict the probability
of a cancer patient's progression-free survival beyond a certain time
horizon [\citet{hari}]. Economists give the probability forecasts
of an
economic rebound or a recession by the end of a fiscal year. Banks are
required by regulators assessing their capital requirements to predict
periodically the risk of default of the loans they make. Engineers are
routinely called upon to predict the survival probability of a system
or infrastructure beyond five or ten years; this includes bridges,
sewer systems and other structures. Finally, lawyers also assess the
probability of particular trial outcome [\citet{fox}] in order to
determine whether to go to trial or settle out of court. This list
would not be complete without mentioning the field that is most
advanced in its daily probability predictions, namely meteorology. In
the past 60 years, remarkable advances in forecasting precipitation
probabilities, temperatures, and rainfall amounts have been made in
terms of breadth and accuracy. \citet{murphy} provide an illuminating
history of the US National Weather Service's transition from
nonprobabilistic to probability predictions and its development of
reliability and accuracy measures for these probability forecasts.
Accuracy assessment is difficult to carry out directly because it
requires comparing a forecaster's predicted probabilities with the
actual but unknown probabilities of the events under study.
Reliability is measured using ``scoring rules,'' which are empirical
distance measures between repeated predicted probabilities of an
event, such as having no rain the next day, and indicator variables
that take on the value 1 if the predicted event actually occurs, and 0
otherwise; see \citet{gneit}, \citet{gneit2} and
\citet
{ranjan} for
recent reviews and developments.

To be more specific, a scoring rule for a sequence of $n$ probability
forecasts $\hp_i, i=1,\ldots,n$, is the average score
$n^{-1}\sum_{i=1}^nL(Y_i,\hp_i)$, where $Y_i=1$ or 0 according to
whether the $i$th event $A_i$ actually occurs or not. An example is the
widely used Brier's score $L(y,\hp)=(y-\hp)^2$ [\citet{brier}]. Noting
that the $Y_i$ are related to the actual but unknown probability $p_i$
via $Y_i\sim\operatorname{Bernoulli}(p_i)$, \citet{cox} proposed
to evaluate
how well the $\hp_i$ predict $p_i$ by using the estimates of
$(\beta_1,\beta_2)$ in the regression model
%
%
\begin{equation}\label{11}
\operatorname{logit}(p_i)=\beta_1+\beta_2 \operatorname{logit}(\hp_i)
\end{equation}
and developed a test of the null hypothesis $(\beta_1,\beta_2)=(0,1)$,
which corresponds to perfect prediction. \citet{spiegel} subsequently
proposed a test of the null hypothesis $H_0\dvtx\hp_i=p_i$ for all
$i=1,\ldots,n$, based on a standardized form (under $H_0$) of Brier's
score. A serious
limitation of this approach is the unrealistic
benchmark of perfect prediction to formulate the null hypothesis, so
significant departures from it are expected when $n$ is large, and they
convey little information on how well the $\hp_i$ predict $p_i$.
Another limitation is the implicit assumption that the $\hp_i$ are
independent random variables, which clearly is violated since $\hp_i$
usually involves previous observations and predictions.

\citet{seillier} have developed a hypothesis testing approach that
removes both limitations in testing the validity of a sequence of
probability forecasts. The forecaster is modeled by a probability
measure under which the conditional probability of the occurrence of
$A_i$ given the $\sigma$-field $\g_{i-1}$ generated by the forecaster's
information set prior to the occurrence of the event is~$\pi_i$. In
this model, the forecaster uses $\hp_i=\pi_i$ as the predicted
probability of $A_i$. As pointed out earlier by \citet{dawid},
this model fits neatly into de~Finetti's (\citeyear{definetti})
framework in which ``the coherent subjectivist Bayesian can be shown to
have a joint probability distribution over all conceivably observable
quantities,'' which is represented by the probability measure $\Pi$ in
the present case. To test if $\Pi$ is ``empirically valid'' based on
the observed outcomes $Y_1,\ldots,Y_n$, \citet{seillier} consider
the null hypothesis $H_0$ that ``the sequence of events is generated by
the same joint distribution from which the forecasts are constructed.''
Under this null hypothesis, $\sum_{i=1}^n\xi_i(Y_i-\pi_i), n\geq1$, is
a martingale with respect to the filtration $\{\g_i\}$ when $\xi_i$ is
$\g_{i-1}$-measurable for all $i$. Assuming certain regularity
conditions on $\xi_i$, they apply the martingale central limit theorem
to show that as $n\to\infty$,
%
%
\begin{equation}\label{new12}
\Biggl\{\sum_{i=1}^n\xi_i(Y_i-\pi_i)\Biggr\}\bigg/\Biggl\{\sum_{i=1}^n\xi_i^2\pi_i(1-\pi
_i)\Biggr\}
^{1/2}\Longrightarrow N(0,1)
\end{equation}
under $H_0$, where $\Longrightarrow$ denotes convergence in
distribution. Since $\pi_i=\hp_i$ in this model of a coherent
forecaster, \citet{seillier} have made use of (\ref{new12}) to
construct various tests of $H_0$. One such test, described at the end
of their Section 6, involves another probability forecast $\hp_i'$,
which is ``based on no more information'' to define $\xi_i$, so that a
significantly large value of the test statistic can be used to reject
$H_0$ in favor of the alternative forecasting model or method.

Hypothesis testing has been extended from testing perfect prediction
or empirical validity of a sequence of probability forecasts to
testing equality of the predictive performance of two forecasts; see
\citet{redel} who extended Spiegelhalter's approach mentioned above.
Testing the equality of predictive performance, measured by some loss
function of the predictors and the realized values, of two forecasting
models or methods has attracted much recent interest in the
econometrics literature, which is reviewed in Section \ref{newsec62}.
In this
paper we develop a new approach to statistical inference, which
involves confidence intervals rather than statistical tests of a null
hypothesis asserting empirical validity of a forecasting model or
method, or equal predictive performance for two forecasting models or
methods. The essence of our approach is to evaluate probability
forecasts via the average loss $L_n=n^{-1}\sum_{i=1}^nL(p_i,\hp_i)$,
where $p_i$ is the actual but unknown probability of the occurrence of
$A_i$. When $L$ is linear in $p_i$, $L(Y_i,\hp_i)$ is an unbiased
estimate of $L(p_i, \hp_i)$ since $E(Y_i|\hp_i)=p_i$. We show in
Section \ref{sec2}, where an overview of loss functions and scoring
rules is
also given, that even for $L$ that is nonlinear in $p_i$ there is a
``linear equivalent'' which carries the same information as $L$ for
comparing different forecasts. In Section \ref{sec3} we make use of this
insight to construct inferential procedures, such as confidence
intervals, for the average loss $L_n$ under certain assumptions and
for comparing the average losses of different forecasts.

Note that we have used $E$ to denote expectation with respect to the
actual probability measure $P$, under which $A_i$ occurs with
probability $p_i$ given the previous history represented by the
$\sigma$-field $\g_{i-1}$, and that we have used $\Pi$ to denote the
probability measure assumed by a coherent Bayesian forecaster whose
probability of occurrence of $A_i$ given $\g_{i-1}$ is $\pi_i$.
Because $\pi_i=\hp_i$ for a coherent Bayesian forecaster,
\citet{seillier} are able to use (\ref{new12}) to test the null
hypothesis of empirical validity of $\Pi$ in the sense that
$E_\Pi(Y_i|\g_{i-1})=\hp_i$, where $E_\Pi$ denotes expectation with
respect to the measure $\Pi$. Replacing $\Pi$ by $P$ is much more
ambitious, but it appears impossible to derive the studentized version
of the obvious estimate $\hl_n=n^{-1}\sum_{i=1}^nL(Y_i,\hp_i)$ and its
sampling distribution under $P$ to perform inference on $L_n$. We
address this difficulty in several steps in Section \ref{sec3}. First we
consider in Section \ref{sec31} the case in which $L(p,\hp)$ is linear
in $p$
and make use of the martingale central limit theorem to prove an
analog of (\ref{new12}) with $p_i$ in place of $\pi_i$ and
$\xi_i=L(1,\hp_i)-L(0,\hp_i)$. Whereas $\pi_i=\hp_i$ under $\Pi$, the
$p_i$ associated with $P$ are unknown parameters that need to be
estimated. Postponing their estimation to Section \ref{sec33}, we
first use
the simple bound $p_i(1-p_i)\leq1/4$ to obtain confidence intervals
for $L_n$ by making use of this analog of (\ref{new12}). In
Section \ref{sec32} we consider the problem of comparing two probability
forecasts via the difference of their average losses, and make use of
the idea of linear equivalents introduced in Section \ref{sec2} to
remove the
assumption of $L(p,\hp)$ being linear in $p$ when we consider
$\Delta_n=n^{-1}\sum_{i=1}^n\{L(p_i,\hp_i')-L(p_i,\hp_i'')\}$. A
variant of $\Delta_n$, called Winkler's skill score in weather
forecasting, is considered in Section \ref{new33}. In Section \ref
{sec33}, we return
to the problem of estimating $p_i(1-p_i)$. Motivated by applications
in which the forecasts are grouped into ``risk buckets'' within which
the $p_i$ can be regarded as equal, Section \ref{sec33} provides two main
results on this problem. The first is Theorem \ref{lem1}, which gives
consistent estimates of the asymptotic variance of $\hdel_n$, or of
$\hl_n$ when $L(p,\hp)$ is linear in $p$, in the presence of risk
buckets with each bucket of size 2 or more. The second, given in
Theorem \ref{lem2}, shows that in this bucket model it is possible to adjust
the Brier score to obtain a consistent and asymptotically normal
estimate of the average squared error loss
$L_n=n^{-1}\sum_{i=1}^n(p_i-\hp_i)^2$. Theorem \ref{lem2} also
provides a
consistent estimate of the asymptotic variance of the adjusted Brier
score when the bucket size is at least 3. In Section \ref{sec35} we
develop an
analog of Theorem \ref{lem1} for the more general setting of
``quasi-buckets,'' for which the $p_i$ within each bin (quasi-bucket)
need not be equal. These quasi-buckets arise in ``reliability
diagrams'' in the meteorology literature. Theorem \ref{thm5} shows
that the
confidence intervals obtained under an assumed bucket model are still
valid but tend to be conservative if the buckets are actually
quasi-buckets. The proofs of Theorems \ref{lem2} and \ref{thm5} are
given in
Section~\ref{appendix}.

Section \ref{sec4} gives a simulation study of the performance of the proposed
methodology, and some concluding remarks and discussion are given in
Section \ref{sec6}. In Section \ref{newsec6} we extend the $Y_i$ from
the case of
indicator variables of events to more general random variables by
modifying the arguments in Section \ref{appendix}, and also show how the
methods and results in Sections \ref{sec32} and \ref{sec33} can be
used to
address related problems in the econometrics literature on the
expected difference in scores between two forecasts, after a brief
review of that literature that has become a major strand of research
in economic forecasts.

\section{Scoring rules and associated loss functions}\label{sec2}

Instead of defining a scoring rule via $L$ (which associates better
forecasts with smaller values of $L$), \citet{gneit} and others assign
higher scores to better forecasts; this is tantamount to using $-L$
instead of $L$ in defining a scoring rule. More generally, considering
$p$ and its forecast $\hp$ as probability measures, they call a
scoring rule $S$ \textit{proper} relative to a class $\mathcal{P}$ of
probability measures if $E_pS(Z,p)\geq E_pS(Z,\hp)$ for all $p$ and
$\hp$ belonging to $\mathcal{P}$, where $Z$ is an observed random
vector (generated from $p$) on which scoring is based. For the
development in the subsequent sections, we find it more convenient to
work with $L$ instead of $-L$ and restrict to $Z=(Y_1,\ldots, Y_n)$ so
that $S(Z,(\hp_1,\ldots,\hp_n))=-n^{-1}\sum_{i=1}^nL(Y_i,\hp_i)$.

The function $L$ in the scoring rule $n^{-1}\sum_{i=1}^nL(Y_i,\hp_i)$
measures the closeness of the probability forecast $\hp_i$ of event
$i$ before the indicator variable $Y_i$ of the event is observed. We
can also use $L$ as a loss function in measuring the accuracy of
$\hp_i$ as an estimate of the probability $p_i$ of event $i$. Besides
the squared error loss $L(p,\hp)=(p-\hp)^2$ used in Brier's score,
another widely used loss function is the Kullback--Leibler divergence,
%
%
\begin{equation}\label{21}
L(p,\hp)=p\log(p/\hp)+(1-p)\log[(1-p)/(1-\hp)],
\end{equation}
which is closely related to the log score introduced by \citet{good},
as shown below. More general loss functions of this type are the
Bregman divergences; see Section 3.5.4 of \citet{grunwald} and Section
2.2 of \citet{gneit}.

We call a loss function $\till(p,\hp)$ a \textit{linear equivalent} of
the loss function $L(p,\hp)$ if $\till(p,\hp)$ is a linear function of
$p$ and
%
%
\begin{equation}\label{22}
L(p,\hp)-\till(p,\hp)\qquad\mbox{does not depend on }\hp.
\end{equation}
For example, $\till(p,\hp)=-2p\hp+\hp^2$ is a linear equivalent of the
squared error loss $(p-\hp)^2$ used by Brier's score. A linear
equivalent $\till$ of the Kullback--Leibler divergence (\ref{21}) is
given by $-\till(p,\hp)=p\log(\hp)+(1-p)\log(1-\hp)$. This is the
conditional expected value (given $\hp$) of
$Y\log(\hp)+(1-Y)\log(1-\hp)$, which is Good's log score. Since the
probability $\hp_i$ is determined before the Bernoulli random variable
$Y_i$ is observed,
%
%
\begin{equation}\label{23}
E\{L(Y_i,\hp_i)|\hp_i,p_i\}=p_iL(1,\hp_i)+(1-p_i)L(0,\hp_i).
\end{equation}
Therefore the conditional expected loss of a scoring rule $L(Y,\hp)$
yields a loss function
%
%
\begin{equation}\label{new24}
\till(p,\hp)=\{L(1,\hp)-L(0,\hp)\}p+L(0,\hp)
\end{equation}
that is linear in $p$. For example, the absolute value scoring rule
$L(Y,\hp)=|Y-\hp|$ is associated with $\till(p,\hp)=p(1-\hp )+(1-p)\hp$
that is linear in each\vspace*{1pt} argument. Using the notation
(\ref{new24}), the scoring rule $L(Y,\hp)$ is proper if
$\till(p,p)\leq\till(p,\hp)$ for all $p, \hp\in[0,1]$, and is strictly
proper if $\min_{0\leq\hp\leq1}\till(p,\hp)$ is uniquely attained at
$p=\hp $. The scoring rule $|Y-\hp|$, therefore, is not proper;
moreover, $|p-\hp|$ does not have a linear equivalent.

\section{A new approach to evaluation of probability forecasts}\label{sec3}

In this section we first consider the evaluation of a sequence of
probability forecasts $\hp_1,\ldots,\hp_n$ based on the corresponding
sequence of indicator variables\break $Y_1,\ldots,Y_n$ that denote whether
the events actually occur or not. Whereas the traditional approach to
evaluating $\hbp=(\hp_1,\ldots,\hp_n)$ uses the scoring rule\break
$n^{-1}\sum_{i=1}^nL(Y_i,\hp_i)$, we propose to evaluate $\hbp$ via
%
%
\begin{equation}\label{31}
L_n=n^{-1}\sum_{i=1}^nL(p_i,\hp_i),
\end{equation}
where $L$ is a loss function, and $p_i$ is the actual probability of
the occurrence of the $i$th event. Allowing the actual probabilities
$p_i$ to be generated by a stochastic system and the forecast $\hp_k$
to depend on an information set $\g_{k-1}$ that consists of the event
and forecast histories and other covariates before $Y_k$ is observed,
the conditional distribution of $Y_i$ given $\g_{i-1}$ and $p_i$ is
Bernoulli($p_i$), and therefore
%
%
\begin{equation}\label{32}
P(Y_i=1|\g_{i-1},p_i)=p_i.
\end{equation}

\subsection{Linear case}\label{sec31}

In view of (\ref{32}), an obvious estimate of the unknown $p_i$
is~$Y_i$. Suppose $L(p,\hp)$ is linear in $p$, as in the case of linear
equivalents of general loss functions. Combining this linearity
property with (\ref{32}) yields
%
%
\begin{equation}\label{33}
E\{L(Y_i,\hp_i)|\g_{i-1},p_i\}=L(p_i,\hp_i),
\end{equation}
and therefore $L(Y_i,\hp_i)-L(p_i,\hp_i)$ is a martingale difference
sequence with respect to $\{\F_i\}$, where $\F_{i-1}$ is the
$\sigma$-field generated by $\g_{i-1}$ and $p_1,\ldots,p_i$. Let
$d_i=L(Y_i,\hp_i)-L(p_i,\hp_i)$. Since $L(y,\hp)$ is linear in $y$, we
can write $L(y,\hp)=a(\hp)y+b(\hp)$. Setting $y=0$ and $y=1$ in this
equation yields $a(\hp)=L(1,\hp)-L(0,\hp)$. Moreover,
$d_i=a(\hp_i)(Y_i-p_i)$. Since $Y_i|\F_i\sim$ Bernoulli$(p_i)$ and
$\hp_i$ is $\F_{i-1}$-measurable,
%
%
\begin{equation}\label{a1}
E(d_i^2|\F_{i-1})=a^2(\hp_i)p_i(1-p_i).
\end{equation}
By (\ref{a1}),
$\sum_1^nE(d_i^2|\F_{i-1})=\sum_1^n\{L(1,\hp_i)-L(0,\hp_i)\}
^2p_i(1-p_i)=O(n)$
a.s., and therefore $n^{-1}\sum_{i=1}^nd_i\to0$ a.s. by the
martingale strong law\break [Williams (\citeyear{williams}),
Section 12.14] proving
$\hl_n-L_n\to0$ a.s. Moreover, if\break
$n^{-1}\sum_1^nE(d_i^2|\F_{i-1})=\sigma_n^2$ converges in probability
to a nonrandom positive constant, then $\sqrt{n}(\hl_n-L_n)/\sigma_n$
has a limiting standard normal distribution by Theorem 1 of
\citet{seillier}. Summarizing, we have the following.
\begin{theorem}\label{thm1}
Suppose $L(p,\hp)$ is linear in $p$. Let
$\hl_n=n^{-1}\sum_{i=1}^nL(Y_i,\hp_i)$, and define $L_n$ by
(\ref{31}). Letting
%
%
\begin{equation}\label{34}
\sigma_n^2=n^{-1}\sum_{i=1}^n\{L(1,\hp_i)-L(0,\hp_i)\}^2p_i(1-p_i),
\end{equation}
assume that $\sigma_n^2=O(1)$ with probability 1. Then $\hl_n-L_n$
converges to 0 with probability 1. If $\sigma_n^2$ converges in
probability to some nonrandom positive constant, then
$\sqrt{n}(\hl_n-L_n)/\sigma_n$ has a limiting standard normal
distribution.
\end{theorem}

To apply Theorem \ref{thm1} to statistical inference on $L_n$, one
needs to
address the issue that $\sigma_n^2$ involves the unknown $p_i$. As
noted in the third paragraph of Section~\ref{sec1}, \citet
{seillier} have
addressed this issue by using $p_i=E_\Pi(Y_i|\g_{i-1})$ under the null
hypothesis $H_0$ that assumes the sequence of events are generated by
the probability measure $\Pi$. This approach is related to the earlier
work of \citet{dawid}, who assumes a ``subjective probability
distribution'' $\Pi$ for the events so that Bayesian forecasts are
given by $\hp_i=\pi_i=E_\Pi(Y_i|\g_{i-1})$. Letting $\xi_t=1$ or 0
according to whether time $t$ is included in the ``test set'' to
evaluate forecasts, he calls the test set ``admissible'' if $\xi_t$
depends only on $\g_{t-1}$, and uses martingale theory to show that
%
%
\begin{equation}\label{62}\quad
\Biggl(\sum_{i=1}^n\xi_iY_i-\sum_{i=1}^n\xi_i\hp_i\Biggr)\bigg/\sum_{i=1}^n\xi_i
\longrightarrow0\qquad\mbox{a.s. }[\Pi]\mbox{ on
}\Biggl\{\sum_{i=1}^n\xi_i=\infty \Biggr\}.
\end{equation}
From (\ref{62}), it follows that for any $0<x<1$, the long-run average
of $Y_i$ (under the subjective probability measure)
associated\vspace*{1pt} with $\hp_i=x$ (i.e., $\xi_i=I_{\{\hp_i=x\}}$)
is equal to $x$ provided that $\sum_{i=1}^nI_{\{\hp_i=x\}}\to\infty$.
Note that Dawid's well-calibration theorem (\ref{62}) involves the
subjective probability measure $\Pi$. \citet{degroot} have noted
that well-calibrated forecasts need not reflect the forecaster's
``honest subjective probabilities,'' that is, need not satisfy Dawid's
coherence criterion \mbox{$\hp_i=\pi_i$}. They therefore use a
criterion called ``refinement'' to compare well-calibrated forecasts.

In this paper we apply Theorem \ref{thm1} to construct confidence intervals
for $L_n$, under the actual probability measure $P$ that generates the
unknown $p_i$ in (\ref{34}). Whereas substituting $p_i$ by $Y_i$ in
$L(p_i,\hp_i)$ leads to a consistent estimate of $L_n$ when $L$ is
linear, such substitution gives 0 as an overly optimistic estimate of
$p_i(1-p_i)=\operatorname{Var}(Y_i|\F_{i-1})$. A conservative confidence
interval for $L_n$ can be obtained by replacing $p_i(1-p_i)$ in
(\ref{34}) by its upper bound $1/4$. In Section \ref{sec33}, we consider
estimation of $\sigma_n^2$ and of $n^{-1}\sum_{i=1}^nL(p_i,\hp_i)$
when $L$ is nonlinear in $p_i$, under additional assumptions on how
the $p_i$ are generated.

\subsection{Application to comparison of probability forecasts}\label{sec32}

Consider two sequences of probability forecasts
$\hbp'=(\hp_1',\ldots,\hp_n')$ and $\hbp''=(\hp_1'',\ldots,\hp
_n'')$ of
$\mathbf{p}=(p_1,\ldots,p_n)$. Suppose\vspace*{1pt} a loss function $L(p,q)$ is
used to
evaluate each forecast, and let $\till(p,q)$ be its linear equivalent.
Since $L(p,q)-\till(p,q)$ does not depend on $q$ in view of (\ref{22}),
it is a function only of $p$, which we denote by $d(p)$. Hence
\begin{eqnarray*}
L(p_i,\hp_i')-L(p_i,\hp_i'')&=&\{\till(p_i,\hp_i')+d(p_i)\}
-\{\till(p_i,\hp_i'')+d(p_i)\}\\
&=&\till(p_i,\hp_i')-\till(p_i,\hp_i'')
\end{eqnarray*}
is a linear function of $p_i$, and therefore we can estimate
$\Delta_n=n^{-1}\sum_{i=1}^n\{L(p_i,\break\hp_i')-L(p_i,\hp_i'')\}$ by the
difference
$n^{-1}\sum_{i=1}^nL(Y_i,\hp_i)-n^{-1}\sum_{i=1}^nL(Y_i,\hp_i')$ of
scores of the two forecasts. Application of Theorem \ref{thm1} then
yields the
following theorem, whose part (ii) is related to (\ref{new24}).
\begin{theorem}\label{thm2}
Let
$\hat\Delta_n=n^{-1}\sum_{i=1}^n\{L(Y_i,\hp_i')-L(Y_i,\hp_i'')\}$
and
%
%
\begin{eqnarray}\label{35}
\delta_i&=&\{L(1,\hp_i')-L(0,\hp_i')\}-\{L(1,\hp_i'')-L(0,\hp_i'')\}
,\nonumber\\[-8pt]\\[-8pt]
s_n^2&=&n^{-1}\sum_{i=1}^n\delta_i^2p_i(1-p_i).
\nonumber
\end{eqnarray}

\mbox{}\hphantom{\textup{i}}\textup{(i)}
Suppose $L$ has a linear equivalent. Letting
$\Delta_n=n^{-1}\sum_{i=1}^n\{L(p_i,\hp_i')-L(p_i,\hp_i'')\}$,
assume that $s_n^2=O(1)$ with probability 1. Then
$\hat\Delta_n-\Delta_n$ converges to 0 with probability 1. If
furthermore $s_n$ converges in probability to some nonrandom
positive constant, then $\sqrt{n}(\hat\Delta_n-\Delta_n)/s_n$ has a
limiting standard normal distribution.

\textup{(ii)} Without assuming that $L$ has a linear equivalent, the
same conclusion as in \textup{(i)} still holds with
$\Delta_n=n^{-1}\sum_{i=1}^n\{\delta_ip_i+L(0,\hp_i')-L(0,\hp _i'')\}$.
\end{theorem}

\subsection{Illustrative applications and skill scores}\label{new33}

As an illustration of Theorem~\ref{thm2}, we compare the Brier scores
$B_k$ for the $k$-day ahead forecasts $\hp_t^{(k)}, 1\leq k\leq7$, for
Queens, NY, provided by US National Weather Service from June 8, 2007,
to March 31, 2009. Table \ref{tableBS} gives the values of $B_1$ and
$B_k-B_{k-1}$ for $2\leq k\leq6$. Using 1/4 to replace $p_i(1-p_i)$ in
(\ref{35}), we can use Theorem \ref{thm2}(i) to construct conservative
95\% confidence intervals for
\[
\Delta(k)=n^{-1}\Biggl\{\sum_{t=1}^n\bigl(p_t-\hp{}^{(k)}_t\bigr)^2-\sum
_{t=1}^n\bigl(p_t-\hp _t^{(k-1)}\bigr)^2\Biggr\},
\]
in which $p_t$ is the actual probability of precipitation on day $t$.
These confidence intervals, which are centered at $B_k-B_{k-1}$, are
given in Table \ref{tableBS}. The results show significant
improvements, by
shortening the lead time by one day, in forecasting precipitation
$k=2,3,4,6$.

%
%
\begin{table}
\caption{Brier scores $B_{1}$ and 95\% confidence intervals for
$\Delta(k)$}\label{tableBS}
\begin{tabular*}{\tablewidth}{@{\extracolsep{\fill}}l d{2.3} d{2.3}
d{2.3} d{2.3} d{2.3} d{2.3}@{}}
\hline
$\bolds{B_{1}}$ & \multicolumn{1}{c}{$\bolds{\Delta(2)}$}
& \multicolumn{1}{c}{$\bolds{\Delta(3)}$}
& \multicolumn{1}{c}{$\bolds{\Delta(4)}$}
& \multicolumn{1}{c}{$\bolds{\Delta(5)}$}
& \multicolumn{1}{c}{$\bolds{\Delta(6)}$}
& \multicolumn{1}{c@{}}{$\bolds{\Delta(7)}$}\\\hline
0.125 &0.021 & 0.012 & 0.020 & 0.010 & 0.015 & 0.007 \\
& \pm0.010 & \pm0.011 & \pm0.012 & \pm0.011 & \pm0.011 &
\pm 0.010\\
\hline
\end{tabular*}
\end{table}

For another application of Theorem \ref{thm2}, we consider
Winkler's (\citeyear{winkler}) skill score. To evaluate weather
forecasts, a skill score that is commonly used is the percentage
improvement in average score over that provided by climatology,
denoted by $\hp_i^c$ and considered as an ``unskilled'' forecaster,
that is,
%
%
\begin{equation}\label{51}
S_n=\Biggl\{n^{-1}\sum_{i=1}^nL(Y_i,\hp_i^c)-n^{-1}\sum_{i=1}^nL(Y_i,\hp
_i)\Biggr\}
\Big/n^{-1}\sum_{i=1}^nL(Y_i,\hp_i^c).
\end{equation}
Climatology refers to the historic relative frequency, also called the
base rate, of precipitation; we can take it to be
$\hp_i^c=(M+1)^{-1}\sum_{t=-M}^0Y_t$. Noting that (\ref{51}) is not a
proper score although it is intuitively appealing, \citet{winkler}
proposed to replace the average climatology score in the denominator
of (\ref{51}) by individual weights $l(\hp_i,\hp_i^c)$, that is,
%
%
\begin{equation}\label{52}
W_n=n^{-1}\sum_{i=1}^n\{L(Y_i,\hp_i)-L(Y_i,\hp_i^c)\}
/l(\hp_i,\hp_i^c),
\end{equation}
where $l(p,c)=\{L(1,p)-L(1,c)\}I_{\{p\geq
c\}}+\{L(0,p)-L(0,c)\}I_{\{p<c\}}$. Theorem~\ref{thm2}(i) can be readily
extended to show that Winkler's score $W_n$ is a consistent estimate
of
%
%
\begin{equation}\label{new53}
w_n=n^{-1}\sum_{i=1}^n\{L(p_i,\hp_i)-L(p_i,\hp_i^c)\}/l(\hp_i,\hp_i^c)
\end{equation}
and that $\sqrt{n}(W_n-w_n)/\tilde{s}_n$ has a limiting standard normal
distribution, where
%
%
\begin{equation}\label{53}
\tilde{s}_n^2=n^{-1}\sum_{i=1}^n\delta_i^2p_i(1-p_i)/l^2(\hp_i,\hp_i^c).
\end{equation}

\begin{figure}

\includegraphics{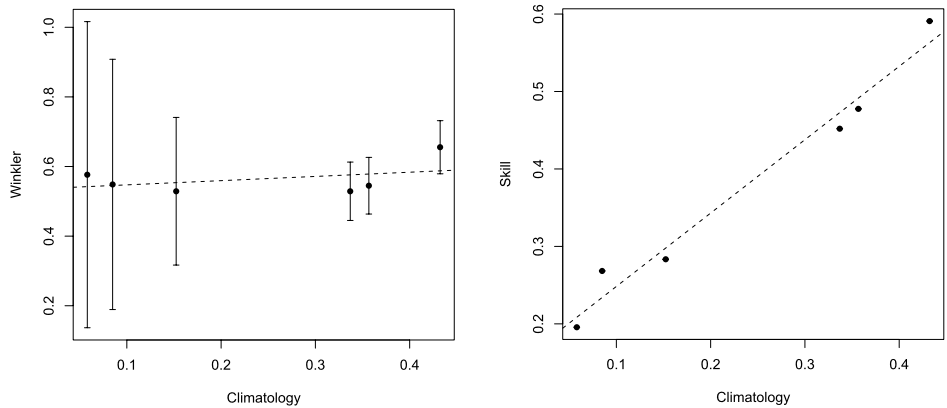}

\caption{The Winkler and skill scores versus climatology.}
\label{doublefig}
\end{figure}

\citet{winkler} used the score (\ref{52}), in which
$L(p,\hp)=(p-\hp)^2$, to evaluate precipitation probability forecasts,
with a 12- to 24-hour lead time, given by the US National Weather
Service for 20 cities in the period between April 1966 and September
1983. Besides the score (\ref{52}), he also computed the Brier score
and the skill score (\ref{51}) of these forecasts and found that both
the Brier and skill scores have high correlations (0.87 and 0.76)
whereas (\ref{52}) has a much lower correlation 0.44 with average
climatology, suggesting that (\ref{52}) provides a better reflection
of the ``skill'' of the forecasts over an unskilled forecasting rule
(based on historic relative frequency). Instead of using correlation
coefficients, we performed a more detailed analysis of Winkler's and
skill scores to evaluate the one-day ahead probability forecasts of
precipitation for six cities: Las Vegas, NV; Phoenix, AZ; Albuquerque,
NM; Queens, NY; Boston, MA; and Portland, OR (listed in increasing
order of relative frequency of precipitation), during the period January
1, 2005, to December 31, 2009. The period January 1, 2002, to December
31, 2004, is used to obtain the past three years' climatology, which
is used as the reference unskilled score in the calculation of the
skill score and Winkler's score (\ref{52}). The left panel of
Figure \ref{doublefig} plots Winkler's score against the relative
precipitation frequency taken from the period January 1, 2005, to
%
%
December 31, 2009, which is simply the percentage of days with rain
during that period and represents the climatology in~(\ref{51}). The
dashed line in the right panel of Figure \ref{doublefig} represents linear
regression of the skill scores (\ref{51}) on climatology and has a
markedly positive slope of 0.95. In contrast, the regression line of
Winkler's scores on climatology, shown in the left panel of
Figure \ref{doublefig}, is relatively flat and has slope 0.12. Unlike skill
scores, Winkler's scores are proper and provide consistent estimates
of the average loss (\ref{new53}) involving the actual daily
precipitation probabilities $p_i$ for each city during the evaluation
period. The vertical bar centered at the dot (representing Winkler's
score) for each city is a 95\% confidence interval for (\ref{new53}),
using a conservative estimate of (\ref{53}) that replaces $p_i(1-p_i)$
by 1/4. The confidence intervals are considerably longer for cities
whose relative frequencies $\hp_i^c$ of precipitation fall below 0.1
because $\delta_i^2/l^2(\hp_i,\hp_i^c)$ tends to be substantially
larger when $\hp_i^c$ is small.

\subsection{Risk buckets and quadratic loss functions}\label{sec33}

Both (\ref{34}) and (\ref{35}) involve $p_i(1-p_i)$, which is the
variance of the Bernoulli random variable $Y_i$. It is not possible to
estimate this variance based on a single observation unless there is
some statistical structure on the $p_i$ to make (\ref{34}) or
(\ref{35}) estimable, and a conservative approach in the absence of
such structure is to use the upper bound $1/4$ for $p_i(1-p_i)$ in
(\ref{34}) or (\ref{35}), as noted in Section \ref{sec31}. One such structure
is that the $p_i$ can be grouped into buckets within which they have
the same value, as in risk assessment of a bank's retail loans (e.g.,
mortgages, automobile loans and personal loans), for which the
obligors are grouped into risk buckets within which they can be
regarded as having the same risk (or more precisely, the same
probability of default on their loans). According to the
Basel Committee on Banking
Supervision [(\citeyear{basel}), page~91] each bank has to use at
least seven risk buckets
for borrowers who have not defaulted and at least one for those who
have defaulted previously at the time of loan application.

A bucket model for risk assessment involves multivariate forecasts for
events $k, 1\leq k\leq K_t$, at a given time $t$. Thus, identifying
the index $i$ with $(t,k)$, one has a vector of probability forecasts
$(\hp_{t,1},\ldots,\hp_{t,K_t})$ at time $t-1$ for the occurrences of
$K_t$ events at time $t$; $K_t=0$ if no forecast is made at time
$t-1$. The information set can then be expressed as $\g_{t-1}$ that
consists of event and forecast histories and other covariates up to
time $t-1$, and therefore conditional on $\g_{t-1}$ and
$p_{t,1},\ldots,p_{t,K_t}$, the events at time $t$ can be regarded as
the outcomes of $K_t$ independent Bernoulli trials with respective
probabilities $p_{t,1},\ldots,p_{t,K_t}$. The bucket model assumes
that, conditional on $\g_{t-1}$ and $p_{t,1},\ldots,p_{t,K_t}$, events
in the same bucket at time $t$ have the same probability of
occurrence. That is, the $p_{t,k}$ are equal for all $k$ belonging to
the same bucket. Let $J_t$ be the number of buckets at time $t$ and
$n_{j,t}$ be the size of the $j$th bucket, $1\leq j\leq J_t$, so that
$n=\sum_{t=1}^T\sum_{j=1}^{J_t}n_{j,t}$. Then\vspace*{1pt} the common $p_{t,k}$ of
the $j$th bucket at time $t$, denoted by $p_t(j)$, can be estimated by
the relative frequency $\bary_t(j)=n_{j,t}^{-1}\sum_{i\in
I_{j,t}}Y_i$, where $I_{j,t}$ denotes the index set for the bucket.
This in turn yields an unbiased estimate
%
%
\begin{equation}\label{36}
\hv_t(j)=n_{j,t}\bary_t(j)\bigl(1-\bary_t(j)\bigr)/(n_{j,t}-1)
\end{equation}
of $p_i(1-p_i)$ for $i\in I_{j,t}$, and we can replace $p_i(1-p_i)$ in
(\ref{34}) or (\ref{35}) by $\hv_t(j)$ for $i\in I_{j,t}$ so that the
results of Theorems \ref{thm1} or \ref{thm2} still hold with these
estimates of
the asymptotic variance, as shown in the following.
\begin{theorem}\label{lem1}
Using the same notation as in the preceding paragraph, suppose
$n_{j,t}\geq2$ for $1\leq j\leq J_t$ and define $\hv_t(j)$ by
(\ref{36}).

\begin{longlist}
\item
Under the same assumptions as in Theorem \ref{thm1}, define
\[
\hat\sigma_n^2=n^{-1}\sum_{t=1}^T\sum_{j=1}^{J_t}\sum_{i\in
I_{j,t}}\{
L(1,\hp_i)-L(0,\hp_i)\}^2\hv_t(j).
\]
Then $\hat\sigma_n^2-\sigma_n^2$ converges to 0 with probability 1.

\item Under the same assumptions as in Theorem \ref{thm2}, $\hat{s}_n^2-s_n^2$
converges to 0 with probability 1, where
$\hat{s}_n^2=n^{-1}\sum_{t=1}^T\sum_{j=1}^{J_t}\sum_{i\in
I_{j,t}}\delta_i^2\hv_t(j)$.
\end{longlist}
\end{theorem}
\begin{pf}
Let $\F_{t-1}$ be the $\sigma$-field generated by $\g_{t-1}$ and
$p_{s,1},\ldots,p_{s,K_s}$ for $s\leq t$. Note that
$\hv_t(j)=\sum_{i\in I_{j,t}}(Y_i-\bary_t(j))^2/(n_{j,t}-1)$ and
that
%
%
\begin{equation}\label{a2}
E(\hv_t(j)|\F_{t-1})=p_t(j)\bigl(1-p_t(j)\bigr),
\end{equation}
which is the variance of $Y_i$ associated with $I_{j,t}$. Therefore
\[
\sum_{j=1}^{J_t}\bigl\{\hv_t(j)-p_t(j)\bigl(1-p_t(j)\bigr)\bigr\}
\biggl\{\sum_{i\in I_{j,t}}[L(1,\hp_i)-L(0,\hp_i)]^2\biggr\}
\]
is a martingale difference sequence with respect to $\{\F_t\}$. Hence
we can apply the martingale strong law as in the proof of Theorem \ref{thm1}
to show that $\hat\sigma_n^2-\sigma_n^2$ converges a.s., and the same
argument also applies to $\hat{s}_n^2-s_n^2$.
\end{pf}

The preceding proof also shows that for the squared error loss
$L(p,\hp)=(p-\hp)^2$, we can estimate (\ref{31}) in the bucket model
by the \textit{adjusted Brier score}
%
%
\begin{equation}\label{new37}
\hl_n-n^{-1}\sum_{t=1}^T\sum_{j=1}^{J_t}n_{j,t}\hv_t(j),
\end{equation}
since $\hl_n=n^{-1}\sum_{i=1}^nL(Y_i,\hp_i)$ is a consistent estimate
of the linear equivalent $n^{-1}\sum_{i=1}^n
(\hp_i^2-2p_i\hp_i+p_i)$, and
$n^{-1}\sum_{t=1}^T\sum_{j=1}^{J_t}n_{j,t}\hv_t(j)$ is a
consistent estimate of $n^{-1}\sum_{i=1}^np_i(1-p_i)$. Consistency of
an estimate $\hat{l}_n$ of $l_n$ means that $\hat{l}_n-l_n$ converges
to 0 in probability as $n\to\infty$. Moreover, the following theorem
shows that
$\sqrt{n}(\hl_n-n^{-1}\sum_{t=1}^T\sum_{j=1}^{J_t}n_{j,t}\hv_t(j)-L_n)$
has a limiting normal distribution in the bucket model and can be
studentized to give a limiting standard normal distribution. Its proof
is given in Section \ref{appendix}.
\begin{theorem}\label{lem2}
Suppose $n_{j,t}\geq2$ for $1\leq j\leq J_t$. Letting
$L(p,\hp)=(p-\hp)^2$, define $L_n$ by (\ref{31}) and the adjusted
Brier score by (\ref{new37}). Let $v_t(j)=p_t(j)(1-p_t(j))$,
%
%
\begin{eqnarray}\label{new38}\quad
\beta_n^2&=&n^{-1}\sum_{t=1}^T\sum_{j=1}^{J_t}\biggl\{v_t(j)\sum_{i\in
I_{j,t}}
(1-2\hp_i)^2\nonumber\\
&&\hphantom{n^{-1}\sum_{t=1}^T\sum_{j=1}^{J_t}\biggl\{}
{}-2v_t(j)\bigl(1-2p_t(j)\bigr)\sum_{i\in I_{j,t}}(1-2\hp_i)\\
&&\hphantom{n^{-1}\sum_{t=1}^T\sum_{j=1}^{J_t}\biggl\{}
{}+n_{j,t}v_t(j)\bigl(1-4v_t(j)\bigr)+2n_{j,t}v^2_t(j)
/(n_{j,t}-1)\biggr\}.\nonumber
\end{eqnarray}
If $\beta_n$ converges in probability to some nonrandom positive
constant, then
$\sqrt{n}(\hl_n-n^{-1}\sum_{t=1}^T\sum_{j=1}^{J_t}n_{j,t}\hv
_t(j)-L_n)/\beta_n$
has a limiting standard normal distribution. Moreover, if $n_j\geq3$
for all $1\leq j\leq J_t$, then $\hat\beta_n-\beta_n$ converges to 0
with probability~1, where
%
%
\begin{eqnarray}\label{new39}
\hat\beta_n^2&=&\frac1{n}\sum_{t=1}^T\sum_{j=1}^{J_t}\biggl\{\hv
_t(j)\sum
_{i\in I_{j,t}}
(1-2\hp_i)^2\nonumber\\
&&\hphantom{\frac1{n}\sum_{t=1}^T\sum_{j=1}^{J_t}\biggl\{}{} -\frac{2n_{j,t}^2}{(n_{j,t}-1)^3}\biggl[\sum_{i\in I_{j,t}}
(1-2\hp_i)\biggr]\biggl[\sum_{i\in
I_{j,t}}\bigl(Y_i-\bary_t(j)\bigr)^3\biggr]\nonumber\\[-8pt]\\[-8pt]
&&\hphantom{\frac1{n}\sum_{t=1}^T\sum_{j=1}^{J_t}\biggl\{}{} +\frac{4n_{j,t}(n_{j,t}-1)}{(n_{j,t}-2)^2}\nonumber\\
&&\hspace*{58pt}{}\times\sum_{i\in I_{j,t}}
\biggl[\frac1{2(n_{j,t}-1)}\sum_{k\in I_{j,t},k\neq i}(Y_i-Y_k)^2-\hv
_t(j)\biggr]^2\biggr\}.\nonumber
\end{eqnarray}
%
\end{theorem}

\subsection{Quasi-buckets and reliability diagrams}\label{sec35}

When the actual $p_{t,i}$ in a bin with index set $I_{j,t}$ are not the
same for all $i\in I_{j,t}$, we call the bin a ``quasi-bucket.'' These
quasi-buckets are the basic components of \textit{reliability diagrams}
that are widely used as graphical tools to evaluate probability
forecasts. In his description of reliability diagrams, Wilks
[(\citeyear{wilks}), Sections 7.1.2, 7.1.3] notes that reliability, or
calibration, relates the forecast to the average observation, ``for
specific values of (i.e., conditional on) the forecast.'' A widely used
approach to ``verification'' of forecasts in meteorology is to group
the forecasts $\hp_i$ into bins so that ``they are rounded
operationally to a finite set of values,'' denoted by
$\hp(1),\ldots,\hp(J)$. Corresponding to each $\hp(j)$ is a set of
observations $Y_i, i\in I_j$, taking the values 0 and~1, where
$I_j=\{i\dvtx\hp_i=\hp(j)\}$. The reliability diagram plots
$\bary(j)=(\sum_{i\in I_j}Y_i)/n_j$ versus $\hp(j)$, where $n_j$ is the
size of $I_j$; see Figure \ref{fig2} in Section \ref{sec4}. Statistical
inference for reliability diagrams has been developed in the
meteorology literature under the assumption of ``independence and
stationarity,'' that is, that $(\hp_i,Y_i)$ are i.i.d. samples from a
bivariate distribution of forecast and observation; see Wilks
[(\citeyear{wilks}), Section 7.9.3] and \citet{brocker}. Under
this assumption, the index sets $I_j$ define a bucket model and a
$(1-\alpha)$-level confidence interval for the common mean $p(j)$ of
the $Y_i$ for $i\in I_j$ is
%
%
\begin{equation}\label{310}
\bary(j)\pm
z_{1-\alpha/2}\bigl\{\bary(j)\bigl(1-\bary(j)\bigr)/n_j\bigr\}^{1/2},
\end{equation}
where $z_q$ is the $q$th quantile of the standard normal distribution.

The assumption of i.i.d. forecast-observation pairs is clearly
violated in weather forecasting, and this has led to the concern that
the confidence intervals given by (\ref{310}) ``are possibly too
narrow'' [\citet{wilks}, page 331]. The temporal dependence
between the
forecast-observation pairs can be handled by incorporating time as in
Section \ref{sec33}. To be specific, let $\hp_{t,k}, k\leq K_t$, be the
probability forecasts, at time $t-1$, of events in the next period. We
divide the set $\{\hp_{t,k}\dvtx k\leq K_t,1\leq t\leq T\}$ into bins
$B_1,\ldots,B_J$, which are typically disjoint sub-intervals of
[0,1]. Let
%
%
\begin{eqnarray}\label{311}
I_{j,t}&=&\{k\dvtx\hp_{t,k}\in B_j\},\nonumber\\
\bary_t(j)&=&\sum_{i\in I_{j,t}}Y_i/n_{j,t},\\
\bary(j)&=&\sum_{t=1}^T\sum_{i\in I_{j,t}}Y_i/n_j,
\nonumber
\end{eqnarray}
where $n_{j,t}$ is the cardinality of $I_{j,t}$ and
$n_j=\sum_{t=1}^Tn_{j,t}$. Note that $n_{j,t}$ and $\bary_t(j)$ have
already been introduced in Section \ref{sec33} and that $\bary(j)$ is the
average of the observations in the $j$th bin, as in (\ref{310}). In
the absence of any assumption on $p_i$ for $i\in I_{j,t}$, these index
sets define quasi-buckets instead of buckets. We can extend the
arguments of Section \ref{sec33} to the general case that makes no assumptions
on the $p_i$ and thereby derive the statistical properties of
$\bary(j)$ without the restrictive assumption of i.i.d.
$(\hp_i,Y_i)$. With the same notation as in Section \ref{sec33}, note
that the
index sets $I_{j,t}$ defined in (\ref{311}) are $\g_{t-1}$-measurable
since the $\hp_{t,k}$ are $\g_{t-1}$-measurable.

Whereas $\bary_t(j)$ is used to estimate the common value of $p_i$ for
$i\in I_{j,t}$ and $\hv_t(j)$, defined in (\ref{36}), is used to
estimate the common value of $p_i(1-p_i)$ in Section~\ref{sec33}, the
$p_i$ in
quasi-buckets need no longer be equal. Replacing $Y_i$ by $p_i$ in
$\bary(j)$ and taking a weighted average of $\hv_t(j)$ over $t$, we
obtain
%
%
\begin{equation}\label{312}
\barp(j)=\frac{\sum_{t=1}^T\sum_{i\in I_{j,t}}p_i}{n_j},\qquad
\hv(j)=\frac
{\sum_{t=1}^Tn_{j,t}\hv_t(j)}{n_j}.
\end{equation}
Instead of (\ref{310}) that is based on overly strong assumptions, we
propose to use
%
%
\begin{equation}\label{313}
\bary(j)\pm z_{1-\alpha/2}\{\hv(j)/n_j\}^{1/2}
\end{equation}
as a $(1-\alpha)$-level confidence interval for $\barp(j)$. Part (iii)
of the following theorem, whose proof is given in Section
\ref{appendix}, shows that the confidence interval tends to be
conservative. Parts (i) and (ii) modify the estimates in Theorem
\ref{lem1} for $\sigma_n^2$ and $s_n^2$ when the $p_i$ in the assumed
buckets turn out to be unequal.
\begin{theorem}\label{thm5}
With the same notation as in Theorem \ref{lem1}, remove the assumption that
$p_i$ are all equal for $i\in I_{j,t}$ but assume that $I_{j,t}$ is
$\g_{t-1}$-measurable for $1\leq j\leq J_t$.

\begin{longlist}
\item Under the assumptions of Theorem \ref{thm1}, let
%
%
\begin{equation}\label{314}\quad
\tsig_n^2=\frac1{n}
\sum_{t=1}^T\sum_{j=1}^{J_t}\sum_{i\in I_{j,t}}\{L(1,\hp_i)-L(0,\hp
_i)\}
^2\bigl(Y_i-\bary_t(j)\bigr)^2
\frac{n_{j,t}}{n_{j,t}-1}.
\end{equation}
Then $\tsig_n^2\geq\sigma_n^2+o(1)$ a.s. Moreover, if the
$p_i$ are equal for all $i\in I_{j,t}$ and $1\leq j\leq J_t$, then
$\tsig_n^2-\sigma_n^2$ converges to 0 a.s.

\item Under the assumptions of Theorem \ref{thm2}, $\tils_n^2\geq s_n^2+o(1)$
a.s., where
%
%
\begin{equation}\label{315}
\tils_n^2=\frac1{n} \sum_{t=1}^T\sum_{j=1}^{J_t}\sum_{i\in
I_{j,t}}\delta_i^2\bigl(Y_i-\bary_t(j)\bigr)^2
\frac{n_{j,t}}{n_{j,t}-1}.
\end{equation}

\item Suppose $J_t=J$ for all $t=1,\ldots,T$. For $1\leq j\leq J$,
define $\bary(j)$ by (\ref{311}), and $\barp(j)$ and $\hv(j)$ by
(\ref{312}), in which $n_j=\sum_{t=1}^Tn_{j,t}$. Let
%
%
\begin{equation}\label{316}
v(j)=n_j^{-1}\sum_{t=1}^T\sum_{i\in I_{j,t}}p_i(1-p_i),
\end{equation}
and let $n=n_1+\cdots+n_J$ be the total sample size. Suppose $n_j/n$
and $v(j)$ converge in probability to nonrandom positive constants as
$n\to\infty$. Then $(n_j/v(j))^{1/2}\{\bary(j)-\barp(j)\}$ has a
limiting standard normal distribution as $n\to\infty$. Moreover,
$\hv(j)\geq v(j)+o_p(1)$ and equality holds if the $p_i$ are equal for
all $i\in I_{j,t}$.
\end{longlist}
\end{theorem}

Note that the numerator of (\ref{36}) is equal to $\sum_{i\in
I_{j,t}}(Y_i-\bary_t(j))^2$. The estimate (\ref{314}) or (\ref{315})
essentially replaces this sum by a weighted sum, using the weights
associated with $p_i(1-p_i)$ in the sum (\ref{34}) or (\ref{35}) that
defines $\sigma_n^2$ or $s_n^2$. The term $n_{j,t}/(n_{j,t}-1)$ in
(\ref{314}) and (\ref{315}) corresponds to the bias correction factor
in the sample variance (\ref{36}). Theorem \ref{thm5} shows that
(\ref{314})
[or (\ref{315})] is still a consistent estimator of $\sigma_n^2$ (or
$s_n^2$) if the bucket model holds, and that it tends to over-estimate
$\sigma_n^2$ (or $s_n^2$) otherwise, erring only on the conservative
side.

\section{Simulation studies}\label{sec4}

The risk buckets in Section \ref{sec33} and the forecasts are usually
based on
covariates. In this section we consider $T=2$ in the case of discrete
covariates so that there are $J_t$ buckets of various sizes for
$n=\sum_{t=1}^2\sum_{j=1}^{J_t}n_{j,t}=300$ probability forecasts
prior to observing the indicator variables $Y_1,\ldots,Y_n$ of the
events. We use the Brier score and its associated loss function
$L(p,\hp)=(p-\hp)^2$ to evaluate the probability forecasts and study
by simulations the adequacy of the estimates $\hat\beta_n^2$ and
$\hat{s}_n^2$ and their use in the normal approximations. The
simulation study covers four scenarios and involves 1,000 simulation
runs for each scenario. Scenario \ref{scen1} considers the Brier score of a
forecasting rule, while Scenarios \ref{scen2}--\ref{scen4} consider the difference of
Brier scores of two forecasts. The bucket sizes and how the $p_i$ and
$\hp_i$ are generated in each scenario are described as follows.
\begin{scen}\label{scen1}
There are ten buckets of size 15 each for each period. The common
values $p_t(j)$ in the buckets are 0.1, 0.25, 0.3, 0.35, 0.4, 0.5, 0.65,
0.7, 0.75 and 0.8, respectively, for $t=1,2$. The probability forecast
$\hp_{t,k}, 1\leq k\leq150$, made at time $t-1$, uses covariate
information to identify the bucket $j$ associated with the $k$th
event at time $t$ and predicts that it occurs with probability
$\bary_{t-1}(j)$, assuming that 150 indicator variables at time 0 are
also observed so that $\bary_0(j)$ is available.
\end{scen}
\begin{scen}\label{scen2}
For each period, there are nine buckets, three of which have size~2
and two of which have size 5; the other bucket sizes are 24, 30, 35
and~45 (one bucket for each size). The bucket probabilities $p_t(j)$
are i.i.d. random variables generated from Uniform (0,1). The
forecast $\hp_{t,k}$ is the same as that in Scenario \ref{scen1}, and there is
another forecast $\hp_{t,k}'=\bary_{t-1}$ that ignores covariate information.
\end{scen}
\begin{scen}\label{scen3}
For each period, there are five buckets of size 30 each, and
$p_t(j)=-0.1+j/5$ for $j=1,\ldots,5$. The two forecasts are the same
as in Scenario~\ref{scen2}.
\end{scen}
\begin{scen}\label{scen4}
This is the same as Scenario \ref{scen3}, except that $p_i$ is uniformly
distributed on $[(j-1)/5,j/5]$ for $i\in I_{j,t}$, that is, the bucket
assumption is only approximately correct.
\end{scen}

\begin{figure}

\includegraphics{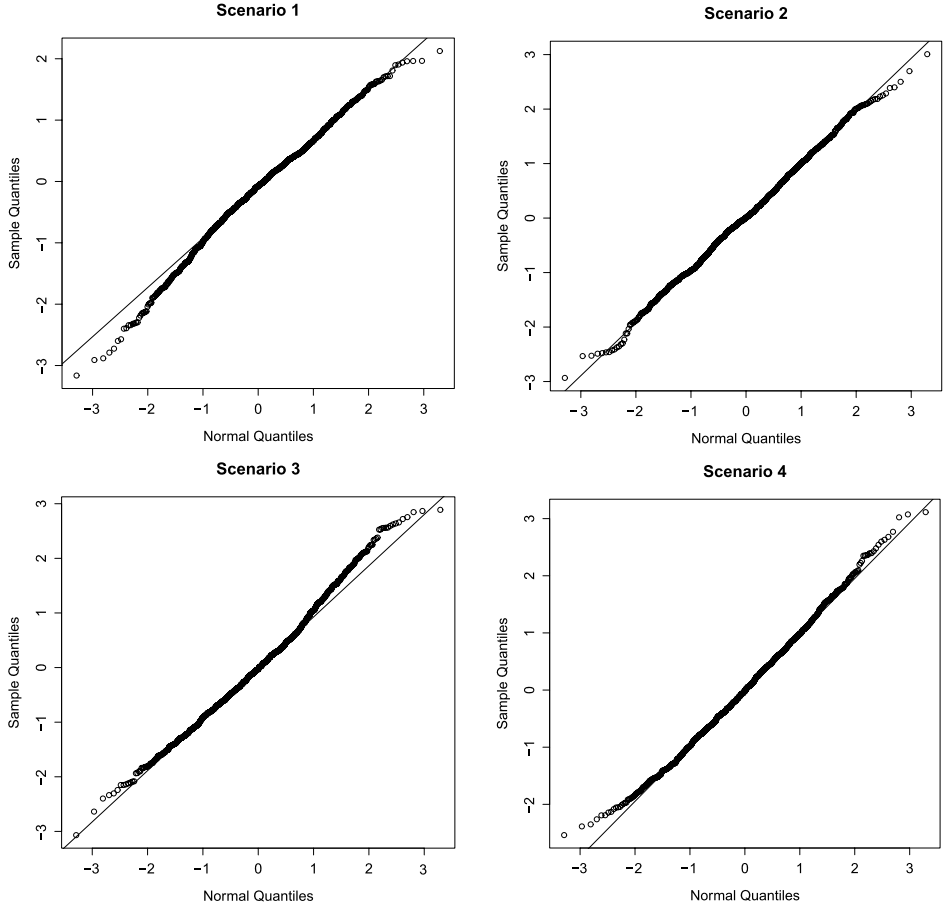}

\caption{Q--Q plots for Scenarios \protect\ref{scen1}--\protect\ref{scen4}.}
\label{MOStable}
\vspace*{-3pt}
\end{figure}

Figure \ref{MOStable} gives the Q--Q plots of
$\sqrt{n}(\hl_n-n^{-1}\sum_{t=1}^2\sum_{j=1}^{J_t}n_{j,t}\hv
_t(j)-L_n)/\hat\beta_n$ for Scenario \ref{scen1} and
$\sqrt{n}(\hat\Delta_n-\Delta_n)/\hat{s}_n$ for Scenarios
\ref{scen2}--\ref{scen4}. Despite the deviation from the assumed bucket
model in Scenario \ref{scen4}, the Q--Q plot does not deviate much from
the $45^\circ$ line. Table \ref{4scenarios} gives the means and
5-number summaries (minimum, maximum, median, 1st and 3rd quartiles) of
$\hat{s}_n/s_n$ for Scenarios \ref{scen2}--\ref{scen4} and
$\hat\beta_n/\beta_n$ for Scenario \ref{scen1}.

%
%

%
%
\begin{table}
\caption{Simulation results for $\hat\beta_n/\beta_n$ (Scenario \protect\ref{scen1})
and $\hat{s}_n/s_n$}\label{4scenarios}
\begin{tabular*}{\tablewidth}{@{\extracolsep{\fill}}lcccccc@{}}
\hline
& \textbf{Min.} & \textbf{1st qrt.} & \textbf{Median} & \textbf{3rd qrt.}
& \textbf{Max.} & \textbf{Mean} \\
\hline
Scenario \ref{scen1}& 0.6397 &1.0840 &1.1810&1.2830 &1.6520&1.1780\\
Scenario \ref{scen2}& 0.7442 &0.9647 &1.0060&1.0490 &1.1970&1.0050\\
Scenario \ref{scen3}& 0.7586 &0.9506 &1.0060&1.0570 &1.2070&1.0010\\
Scenario \ref{scen4}& 0.7420 &0.9661 &1.0180&1.0730 &1.2240&1.0160\\
\hline
\end{tabular*}
\vspace*{-3pt}
\end{table}

\begin{figure}

\includegraphics{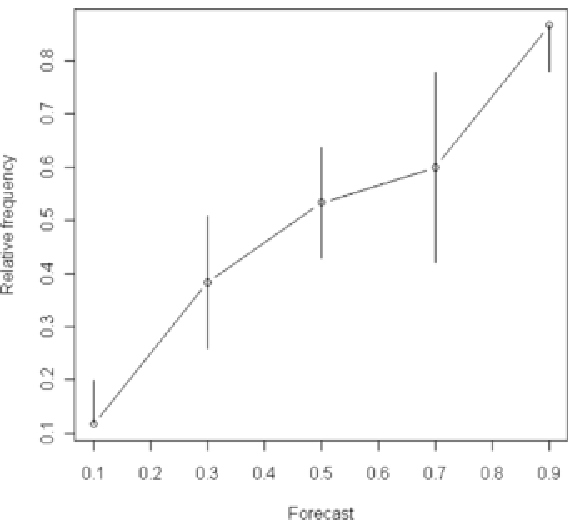}

\caption{Reliability diagram for the forecasts $\hp_{t,k}$. At the
midpoint of each of the five bins $[(j-1)/5,j/5], j=1,\ldots,5$, a
95\% confidence interval, centered at $\bary(j)$, for $\barp(j)$ is
shown; only the upper (or lower) half of the interval is shown at
$j=1$ (or $5$) to keep the range of the vertical axis between $0.1$
and $0.9$.}
\label{fig2}
\vspace*{-12pt}
\end{figure}

To illustrate the reliability diagram and the associated confidence
intervals (\ref{313}) in Section \ref{sec35}, we use one of the
simulated data
sets in Scenario \ref{scen4} to construct the reliability diagram for the
forecasts $\hp_{t,k}$ $(t=1,2;k=1,\ldots,5)$, grouping the forecasts
over time in the bins $[(j-1)/5,j/5], j=1,\ldots,5$, that are natural
for this scenario. The diagram is given in Figure \ref{fig2}.
%
%
Table \ref{table3} gives the means, standard deviations (SD), and
5-number summaries of $\bary(j), \barp(j), \hv(j)$ and $v(j)$ defined
in (\ref{311}), (\ref{312}) and (\ref{316}) based on the 1,000
simulations. In particular, it shows that $\hv(j)$ tends to
%
%
\begin{table}
\caption{Simulation results for $\barp(j),\bary(j),v(j)$ and $\hv(j)$}
\label{table3}
\begin{tabular*}{\tablewidth}{@{\extracolsep{\fill}}lccccccc@{}}
\hline
&\textbf{Min} & \textbf{1st qrt.} & \textbf{Median}
& \textbf{3rd qrt.} & \textbf{Max} & \textbf{Mean} & \textbf{SD}\\
\hline
$\barp(1)$&0.101& 0.101& 0.101& 0.168& 0.234& 0.121& 0.033\\
$\bary(1)$&0.017& 0.083& 0.117& 0.156& 0.350& 0.123& 0.051\\
$v(1)$ &0.087& 0.087& 0.087& 0.127& 0.167& 0.100& 0.020\\
$\hv(1)$ &0.067& 0.089& 0.106& 0.132& 0.233& 0.106& 0.037\\[6pt]
$\barp(2)$&0.101& 0.300& 0.300& 0.355& 0.515& 0.320& 0.049\\
$\bary(2)$&0.050& 0.267& 0.317& 0.378& 0.633& 0.319& 0.089\\
$v(2)$ &0.087& 0.207& 0.207& 0.207& 0.247& 0.209& 0.015\\
$\hv(2)$ &0.048& 0.208& 0.221& 0.239& 0.259& 0.213& 0.034\\[6pt]
$\barp(3)$&0.300& 0.515& 0.515& 0.577& 0.701& 0.527& 0.058\\
$\bary(3)$&0.217& 0.467& 0.533& 0.589& 0.833& 0.529& 0.096\\
$v(3)$ &0.206& 0.233& 0.247& 0.247& 0.247& 0.239& 0.011\\
$\hv(3)$ &0.144& 0.240& 0.249& 0.254& 0.259& 0.244& 0.015\\[6pt]
$\barp(4)$&0.515& 0.659& 0.701& 0.701& 0.906& 0.690& 0.052\\
$\bary(4)$&0.367& 0.633& 0.689& 0.750& 1.000& 0.687& 0.090\\
$v(4)$ &0.082& 0.206& 0.206& 0.206& 0.247& 0.204& 0.021\\
$\hv(4)$ &0.063& 0.207& 0.217& 0.236& 0.259& 0.211& 0.035\\[6pt]
$\barp(5)$&0.769& 0.906& 0.906& 0.906& 0.906& 0.895& 0.026\\
$\bary(5)$&0.733& 0.867& 0.900& 0.933& 1.000& 0.892& 0.049\\
$v(5)$ &0.082& 0.082& 0.082& 0.082& 0.164& 0.088& 0.016\\
$\hv(5)$ &0.077& 0.084& 0.093& 0.120& 0.202& 0.096& 0.037\\
\hline
\end{tabular*}
\vspace*{-12pt}
\end{table}
over-estimate $v(j)$. Moreover, the probability of coverage of the 95\%
interval (\ref{313}) for $\barp(j)$, evaluated by averaging over the
1,000 simulations, is 0.949, 0.947, 0.944, 0.940 and 0.928, for
$j=1,\ldots,5$, respectively, suggesting that the results of Theorem
\ref{thm5} still apply even for moderate sample sizes. We do not
consider the second forecast $\hp_{t,k}'=\bary_{t-1}$ to illustrate
reliability diagrams because by the central limit theorem, the
$\hp_{t,k}'$ are concentrated around 0.5 and nearly all of the
forecasts lie in the bin $[0.4,0.6]$.

\section{\texorpdfstring{Proofs of Theorems \protect\ref{lem2} and \protect\ref
{thm5}}{Proofs of Theorems 4 and 5}}\label{appendix}

Re-labeling the $Y_i$ as $Y_{t,1},\ldots,Y_{t,K_t}$, we note that
conditional on $\F_{t-1}, \{Y_{t,k}\dvtx1\leq k\leq K_t\}$ is a set of
independent Bernoulli random variables with respective parameters
$p_{t,1},\ldots,p_{t,K_t}$. This point, which has been noted in the
second paragraph of Section \ref{sec33} and will be discussed further in
Section \ref{sec6}, explains why we can first derive the result for the
special case in which $Y_{t,k}$ are independent and then modify the
argument by conditioning on $\F_{t-1}$ and appealing to martingale
theory. As an\vadjust{\goodbreak} illustration, note that if $Y_i, i\in I_{j,t}$, are
i.i.d. Bernoulli random variables with common parameter $p_t(j)$,
then $\hv_t(j)$ defined in~(\ref{36}) is an unbiased estimate of
$p_t(j)(1-p_t(j))$ and one can use the classical strong law of large
numbers to derive the result. The proof of Theorem \ref{lem1}
basically shows
that $\hv_t(j)$ is ``conditionally unbiased'' given $\F_{t-1}$ in the
sense of (\ref{a2}) and then uses the martingale strong law to derive
the result. To prove Theorem \ref{lem2}, we extend this idea to obtain a
conditionally unbiased estimate of Var($\hv_t(j)$) by first
considering the i.i.d. case: let $X_1,\ldots,X_m$ be i.i.d. random
variables. As is well known, the sample variance
$\hv=\sum_{i=1}^m(X_i-\bar{X})^2/(m-1)$ is a $U$-statistic of order 2,
with kernel $h(X_i,X_k)=(X_i-X_k)^2/2$. \citet{arvesen} has shown that
an unbiased estimate of the variance of the $U$-statistic is the
jackknife estimate
%
%
\begin{equation}\label{a8}
\frac{4(m-1)}{m(m-2)^2}\sum_{i=1}^m\Biggl\{\frac1{m-1}
\mathop{\sum_{k=1}}_{k\neq i}^mh(X_i,X_k)-\hv\Biggr\}^2.
\end{equation}
\begin{pf*}{Proof of Theorem \ref{lem2}}
Use $L(p,\hp)=(p-\hp)^2$ to express $n\{\hl_n-n^{-1}
\times\break\sum_{t=1}^T\sum_{j=1}^{J_t}n_{j,t}\hv_t(j)-L_n\}$ as
%
%
\begin{equation}\label{a3}
\sum_{t=1}^T\sum_{k=1}^{K_t}(1-2\hp_{t,k})(Y_{t,k}-p_{t,k})-
\sum_{t=1}^T\sum_{j=1}^{J_t}n_{j,t}
\bigl[\hv_t(j)-p_t(j)\bigl(1-p_t(j)\bigr)\bigr],\hspace*{-24pt}
\end{equation}
which is the difference of two martingales and is therefore a
martingale. To compute the conditional variance (or predictable
variation) of (\ref{a3}), we can use the ``angle bracket'' notation
and formulas for predictable variation and covariation
[\citet{williams}, Section 12.12] to obtain
%
%
\begin{eqnarray}\label{a4}
&&\Biggl\langle\sum_{t=1}^T\sum_{k=1}^{K_t}(1-2\hp
_{t,k})(Y_{t,k}-p_{t,k})\Biggr\rangle\nonumber\\
&&\qquad=\sum_{t=1}^T\sum_{k=1}^{K_t}(1-2\hp
_{t,k})^2E\bigl((Y_{t,k}-p_{t,k})^2|\F
_{t-1}\bigr)\\
&&\qquad=\sum_{t=1}^T\sum_{j=1}^{J_t}\biggl(\sum_{i\in I_{j,t}}(1-2\hp
_i)^2\biggr)p_t(j)\bigl(1-p_t(j)\bigr),\nonumber
\\
%
%
\label{a5}
&&\Biggl\langle\sum_{t=1}^T\sum_{k=1}^{K_t}(1-2\hp
_{t,k})(Y_{t,k}-p_{t,k}),\sum
_{t=1}^T\sum_{j=1}^{J_t}n_{j,t}
\bigl[\hv_t(j)-p_t(j)
\bigl(1-p_t(j)\bigr)\bigr]\Biggr\rangle
\hspace*{-25pt}\nonumber\\[-8pt]\\[-8pt]
&&\qquad
=\sum_{t=1}^T\sum_{j=1}^{J_t}\biggl[\sum_{i\in I_{j,t}}(1-2\hp
_i)\biggr]p_t(j)\bigl(1-p_t(j)\bigr)\bigl(1-2p_t(j)\bigr),
\nonumber
\\
%
%
\label{a6}
&&\Biggl\langle\sum_{t=1}^T\sum_{j=1}^{J_t}n_{j,t}\bigl[\hv
_t(j)-p_t(j)\bigl(1-p_t(j)\bigr)\bigr]\Biggr\rangle\nonumber\\
&&\qquad
=\sum_{t=1}^T\sum_{j=1}^{J_t}\bigl\{
n_{j,t}p_t(j)\bigl(1-p_t(j)\bigr)\bigl[1-4p_t(j)\bigl(1-p_t(j)\bigr)\bigr]\\
&&\hspace*{70.5pt}\qquad\quad{}+2n_{j,t}p_t^2(j)\bigl(1-p_t(j)\bigr)^2/(n_{j,t}-1)\bigr\}.\nonumber
\end{eqnarray}
Combining (\ref{a4}), (\ref{a5}) and (\ref{a6}) yields formula
(\ref{new38}) for the conditional variance of (\ref{a3}) divided by $n$.

In view of (\ref{a2}), $\hv_t(j)$ is a conditionally unbiased estimate
of $p_t(j)(1-p_t(j))$ given $\F_{t-1}$. If $Y\sim$ Bernoulli$(p_i)$,
then $E(Y-p)^3=p(1-p)(1-2p)$. Hence a conditionally unbiased estimate
of $p_t(j)(1-p_t(j))(1-2p_t(j))$ given $\F_{t-1}$ is
%
%
\begin{equation}\label{a7}
[n_{j,t}^2/(n_{j,t}-1)^3]\sum_{i\in I_{j,t}}\bigl(Y_i-\bary_t(j)\bigr)^3,
\end{equation}
analogous to (\ref{a2}). Replacing $p_t(j)(1-p_t(j))(1-2p_t(j))$ in
(\ref{a5}) by (\ref{a7}) and multiplying (\ref{a5}) by $-2/n$ gives
the second summand of (\ref{new39}). Note that the first summand of
(\ref{new39}) corresponds to replacing $p_t(j)(1-p_t(j))$ in
(\ref{a4}) by $\hv_t(j)$. The last summand of (\ref{new39})
corresponds to using the jackknife estimate (\ref{a8}) to estimate the
conditional variance of $\hv_t(j)$ given $\F_{t-1}$. Since $\{
Y_i,i\in
I_{j,t}\}$ is a set of i.i.d. random variables conditional on
$\F_{t-1}$, the jackknife estimate is conditionally unbiased given
$\F_{t-1}$; see the paragraph preceding the proof of this theorem. The
rest of the argument is similar to that of Theorem \ref{lem1}.
\end{pf*}
\begin{pf*}{Proof of Theorem \ref{thm5}}
We first prove (iii). Using the notation in the paragraph preceding
the proof of Theorem \ref{lem2}, recall that conditional on $\F
_{t-1}$, the
$Y_{t,k}$ are independent Bernoulli($p_{t,k}$) random variables. Since
$I_{j,t}$ is $\F_{t-1}$-measurable, it follows that $\sum_{i\in
I_{j,t}}(Y_i-p_i)$ is a martingale difference sequence with respect
to $\{\F_t\}$ and $E\{[\sum_{i\in
I_{j,t}}(Y_i-p_i)]^2|\F_{t-1}\}=\sum_{i\in
I_{j,t}}p_i(1-p_i)$. Since $n^{-1}\sum_{t=1}^T\sum_{i\in
I_{j,t}}p_i(1-p_i)$ converges in probability to a nonrandom
positive constant as $n\to\infty$, we can apply the martingale central
limit theorem as in the proof of Theorem \ref{thm1} to conclude that
\[
\frac{\sum_{t=1}^T\sum_{i\in
I_{j,t}}(Y_i-p_i)}{\{\sum_{t=1}^T\sum_{i\in
I_{j,t}}p_i(1-p_i)\}^{1/2}}\Longrightarrow N(0,1)
\]
proving the first part of (iii).

To prove the second part of (iii), and also (i) and (ii), we first
show that for any nonnegative $\F_{t-1}$-measurable random variables
$w_{t,1},\ldots,w_{t,K_t}$,
%
%
\begin{eqnarray}\label{59}
&&E\biggl\{\sum_{k\in
I_{j,t}}w_{t,k}
\bigl(Y_{t,k}-\bary_t(j)\bigr)^2|\F_{t-1}\biggr\}\nonumber\\[-8pt]\\[-8pt]
&&\qquad\geq\sum_{k\in
I_{j,t}}(1-n_{j,t}^{-1})w_{t,k}p_{t,k}(1-p_{t,k}),\nonumber
\end{eqnarray}
in which $\sum_{k\in I_{j,t}}$ means $\sum_{i\in I_{j,t}}$ when $i$ is
represented as $(t,k)$; see the second paragraph of
Section \ref{sec33}. Define $\bary_t(j), \bary(j)$ and $\barp(j)$
as in
(\ref{311}) and (\ref{312}), and let $\barp_t(j)=(\sum_{k\in
I_{j,t}}p_{t,k})/n_{j,t}$. From the decomposition
%
%
\begin{equation}\label{510}
Y_{t,k}-\bary_t(j)=(Y_{t,k}-p_{t,k})+\bigl(p_{t,k}-\barp_t(j)\bigr)+\bigl(\barp
_t(j)-\bary_t(j)\bigr),
\end{equation}
it follows that the left-hand side of (\ref{59}) is equal to
%
%
\begin{eqnarray}\label{511}
&&\sum_{k\in I_{j,t}}w_{t,k}E[(Y_{t,k}-p_{t,k})^2|\F_{t-1}]+\sum_{k\in
I_{j,t}}w_{t,k}\bigl(p_{t,k}-\barp_t(j)\bigr)^2\nonumber\\
&&\qquad{}+\sum_{k\in I_{j,t}}w_{t,k}E\bigl[\bigl(\bary_t(j)-\barp_t(j)\bigr)^2|\F_{t-1}\bigr]\\
&&\qquad{}-2E\biggl[\bigl(\bary_t(j)-\barp_t(j)\bigr)\sum_{k\in
I_{j,t}}w_{t,k}(Y_{t,k}-p_{t,k})\biggr]\nonumber
\end{eqnarray}
by using\vspace*{1pt} the fact that conditional on $\F_{t-1}$ the $Y_{t,k}$ are
independent Bernoulli. Since $\bary_t(j)-\barp_t(j)=\sum_{k\in
I_{t,j}}(Y_{t,k}-p_{t,k})/n_{j,t}$, we can\vspace*{1pt} use this fact again to
combine the last two terms of (\ref{511}) into
%
%
\begin{equation}\label{512}
-\sum_{k\in I_{j,t}}(w_{t,k}/n_{j,t})E[(Y_{t,k}-p_{t,k})^2|\F_{t-1}].
\end{equation}
Since $w_{t,k}\geq0$, we can drop the second term in (\ref{511}) to
obtain (\ref{59}) from (\ref{511}) and (\ref{512}). Moreover, since
this term is actually 0 when the $p_{t,k}$ are all equal for $k\in
I_{j,t}$, equality holds in (\ref{59}) in this case.

Let $w_{t,k}=n_{j,t}/(n_{j,t}-1)$. Then (\ref{59}) reduces to
%
%
\begin{equation}\label{513}
E(n_{j,t}\hv_t(j)|\F_{t-1})\geq\sum_{k\in I_{j,t}}p_{t,k}(1-p_{t,k}).
\end{equation}
Under the assumptions of part (iii) of the theorem, we can apply the
martingale strong law to obtain
%
%
\begin{equation}\label{514}\qquad
\sum_{t=1}^T\{n_{j,t}\hv_t(j)-E(n_{j,t}\hv_t(j)|\F_{t-1})\}
/n_j\longrightarrow0 \qquad\mbox{a.s.
on }\{n_j\to\infty\}.
\end{equation}
Combining (\ref{513}) with (\ref{514}) yields $\hv(j)\geq
v(j)+o_p(1)$, with equality when the $p_{t,k}$ are all equal for $k\in
I_{j,t}$.

To prove part (i) of the theorem, put
$w_{t,k}=\{L(1,\hp_{t,k})-L(0,\hp_{t,k})\}^2n_{j,t}/\break(n_{j,t}-1)$ in
(\ref{59}) and then use the same argument as in the preceding
paragraph. The proof of part (ii) is similar.
\end{pf*}

\section{Extensions and connections to forecast comparison in
econometrics}\label{newsec6}

Our new approach to evaluating probability forecasts in Section \ref
{sec3} is
based on consistent and asymptotically normal estimates of the average
loss\break $n^{-1}\sum_{i=1}^nL(p_i,\hp_i)$, without any assumptions on how
the observed indicator variables $Y_i$ and their forecasts $\hp_i$ are
generated. The key to this approach is that conditional on $\F_{i-1}$,
$Y_i$ is Bernoulli$(p_i)$, and therefore martingale arguments can be
used to derive the results in Section \ref{sec3}. In Section \ref
{newsec61} we show
how this approach can be extended to more general random variables
$Y_i$. As shown in (\ref{23}), when $Y_i$ is an indicator variable,
the conditional expectation of the score $L(Y_i,\hp_i)$ given
$\F_{i-1}$ is a linear function of $p_i$, but this does not extend to
more general random variables $Y_i$. In Section \ref{newsec62} we
review the
recent econometrics literature on testing the equality of the expected
scores of two forecasts and discuss an alternative approach to
statistical inference on the expected difference in average scores of
two forecasts.

\subsection{Extensions to general predictands}\label{newsec61}

A characteristic of $(\hp_i,Y_i)$ in probability forecasting is that
$E(Y_i|\F_{i-1})=p_i$ while the $\g_{i-1}$-measurable forecast $\hp_i$
is an estimate of $p_i$. The theorems in Section \ref{sec3} and their
martingale proofs in Section~\ref{appendix} can be easily extended to general
random variables $Y_i$ when the loss function is of the form
$L(\mu_i,\hmu_i)$, where $\mu_i=E(Y_i|\F_{i-1})$ and $\hmu_i$ is a
forecast of $Y_i$ given~$\F_{i-1}$. Although $Y_i|\F_{i-1}\sim$
Bernoulli$(p_i)$ in Section \ref{sec3}, no parametric assumptions are actually
needed when we use a loss function of the form $L(\mu_i,\hmu_i)$. As
in (\ref{22}), such loss function is said to have a linear equivalent
$\till$ if
%
%
\begin{equation}\label{n61}\qquad
\till(y,\hat{y})\mbox{ is linear in $y$}\quad\mbox{and}\quad L(y,\hat{y})-\till
(y,\hat{y})\mbox{
does not depend on }\hat{y}.
\end{equation}
The bucket model in Section \ref{sec33} can be extended so that
$Y_{t,k}|\F_{t-1}$ have the same mean and variance for all $(t,k)$
belonging to the same bucket. In place of (\ref{36}), we now use
%
%
\begin{equation}\label{n62}
\hv_t(j)=\sum_{i\in I_{j,t}}\bigl(Y_i-\bary_t(j)\bigr)^2/(n_{j,t}-1)
\end{equation}
as an unbiased estimate of the common conditional variance of $Y_i$
given $\F_{i-1}$ for $i=(t,k)\in I_{j,t}$, using the same notation as
that in the proof of Theorem \ref{lem2}. While the extension of Theorem
\ref{lem1}
only needs the first two moments of $Y_{t,k}|\F_{t-1}$ to be equal for
all $(t,k)$ belonging to the same bucket, Theorem \ref{lem2} can also be
extended by assuming the first four moments of $Y_{t,k}|\F_{t-1}$ to
be equal for all $(t,k)$ belonging to the same bucket, by using
Arvesen's (\citeyear{arvesen}) jackknife estimate of the variance of a
$U$-statistic.

Clearly (\ref{510}), (\ref{511}) and (\ref{512}) also hold with
$p_{t,k}$ and $\barp_t(j)$ replaced by $\mu_{t,k}$ and $\barm_t(k)$,
so Theorem \ref{thm5} can likewise be extended to quasi-buckets and
reliability diagrams for the predicted means $\hmu_{t,k}$. For sample
means in the case of independent observations within each bucket, this
extension of Theorem \ref{thm5} can be viewed as a corollary of the analysis
of variance. In fact, the proof of Theorem \ref{thm5} uses martingale
arguments and conditioning to allow dependent observations in each
(quasi-)bucket.

\subsection{Inference on expected difference in average scores of two
forecasts}\label{newsec62}

When the $Y_i$ are indicator variables of events, Theorem \ref{thm2}(ii)
establishes asymptotic normality for the difference
$\hdel_n=n^{-1}\sum_{i=1}^n\{L(Y_i,\hp_i')-L(Y_i,\hp_i'')\}$ in
average scores between two forecasts, from which one can perform
inference on
%
%
\begin{eqnarray}\label{n63}
\Delta_n&=&n^{-1}\sum_{i=1}^nE\{L(Y_i,\hp_i')-L(Y_i,\hp_i'')|\F
_{i-1}\}\nonumber\\[-8pt]\\[-8pt]
&=&n^{-1}\sum_{i=1}^n\{\delta_ip_i+L(0,\hp_i')-L(0,\hp_i'')\},
\nonumber
\end{eqnarray}
where
$\delta_i=\{L(1,\hp_i')-L(0,\hp_i')\}-\{L(1,\hp_i'')-L(0,\hp_i'')\}$.
This simplicity, however, does not extend to general $Y_i$.

Proper scoring rules for probability forecasts of categorical and
continuous variables $Y_i$ have been an active area of research; see
the review by \citet{gneit}. Another active area of research is related
to the extension of $\hdel_n$ to general $Y_i$ in the econometrics
literature, beginning with the seminal paper of \citet{diebold}. They
consider the usual forecast errors $e_t:=Y_t-\hY_{t|t-1}$ in time
series analysis, where $\hY_{t|t-1}$ is the one-step ahead forecast of
$Y_t$ based on observations up to time $t-1$. Unlike a probability
forecast that gives a predictive distribution of $Y_t$ as in
\citet{gneit}, $\hY_{t|t-1}$~is a nonprobabilistic forecast that
predicts the value of $Y_t$ [see \citet{wilks}, Section 7.3]. The score
used by \citet{diebold} is of the form $L(Y_t,\hY
_{t|t-1})=g(e_t)$, and
they consider the \textit{average loss differential}
%
%
\begin{equation}\label{n64}
\hdel_n=n^{-1}\sum_{t=1}^n\{L(Y_t,\hY_{t|t-1}')-L(Y_t,\hY
_{t|t-1}'')\}
\end{equation}
between two forecasts $Y_{t|t-1}'$ and $\hY_{t|t-1}'', 1\leq t\leq
n$. Assuming a probability measure $Q$ under which
$d_t:=g(e_t')-g(e_t'')$ is covariance stationary with absolutely
summable autocovariances $\gamma_k$ so that
$f(0):=\sum_{k=-\infty}^\infty\gamma_k/(2\pi)$ is the spectral density
at frequency 0, they use the asymptotic normality of $\hdel_n$ under
the null hypothesis $H_0\dvtx E_Q(d_t)=0$ and a window estimate $\hf(0)$
of $f(0)$ so that the test statistic
$\sqrt{n}\hdel_n/(2\pi\hf(0))^{1/2}$ has a limiting standard normal
distribution under $H_0$ as $n\to\infty$. This aysmptotic normality
result, however, requires additional assumptions, such as stationary
mixing, which they do not mention explicitly. Their work has attracted
immediate attention and spawned many subsequent developments in the
econometrics literature on this topic.

\citet{giacomini}, hereafter abbreviated as G\&W, review some of the
developments\vspace*{2pt} and propose a refinement of $H_0$ for which the
asymptotic normality of $\hdel_n$ can be established under precisely
stated conditions that can also allow nonstationarity. They formulate
the null hypothesis of equal predictive ability of two forecasting
models or methods as ``a problem of inference about conditional
expectations of forecasts and forecast errors that nests the
unconditional expectations that are the sole focus of the existing
literature.'' The econometrics literature they refer to is primarily
concerned with ``forecast models;'' thus $Q$ in the previous paragraph
is the probability measure associated with the forecast model being
evaluated, or with a more general model than the two competing
forecast models whose predictive abilities are compared.

G\&W evaluate not only the forecasting model but also the forecasting
method, which includes ``the forecasting model along with a number of
choices,'' such as the estimation procedure and the window of past
data, used to produce the forecast. They consider $k$-step ahead
forecasts, for which $\hY_{t|t-1}$ is replaced by $\hY_{t|t-k}$, and
assume that the forecasts are based on finite-memory models involving
unknown parameters, that is,
%
%
\begin{equation}\label{n65}
\hY_{t+k|t}=h(Y_t,\ldots,Y_{t-m+1},\bx_t,\ldots,\bx_{t-m+1};\bbet),
\end{equation}
where $h$ is a known function, $m$ is the order of the model, $\bx_t$
is a covariate vector at time $t$ and $\bbet$ is a parameter vector
to be estimated from some specified window of past data. Their
formulation generalizes that of \citet{west} who considers regression
models. Whereas West assumes that the data are actually generated by
the regression model with true parameter $\bbet^*$, G\&W allow model
misspecification, and therefore their assumptions do not involve
$\bbet^*$. They consider two such nominal models, resulting in the
forecasts $\hY_{t|t-k}'$ and $\hY_{t|t-k}''$ that use the same
covariates but different estimates $\hat\bbet_t'$ and $\hat\bbet_t''$.
Their null hypothesis
%
%
\begin{equation}\label{n66}
H_0\dvtx E\{L(Y_t,\hY_{t|t-k}')-L(Y_t,\hY_{t|t-k}'')|\g_{t-k}\}=0
\qquad\mbox{a.s. }\forall t\geq1
\end{equation}
seems to be stronger than
$E L(Y_t,Y_{t|t-k}')=E L(Y_t,Y_{t|t-k}'')$ $\forall t$ considered by
\citet{diebold} for the case $k=1$. On the other hand, (\ref
{n66}) in
the case $k=1$ just says that $L(Y_t,Y_{t|t-1}')-L(Y_t,Y_{t|t-1}'')$
is a martingale\vspace*{1pt} difference sequence under $H_0$ so that the martingale
central limit theorem can be applied to derive the limiting
$\chi^2$-distribution of G\&W's test statistics under $H_0$. Unlike
\citet{diebold} for the case $k=1$, G\&W do not use test
statistics of
the form (\ref{n64}) and their test statistics involve more
complicated weighted sums of
$L(Y_t,\hY_{t|t-k}')-L(Y_t,\hY_{t|t-k}'')$. These\vspace*{1pt} weights are chosen
to improve the power of the test and require additional mixing and
moment assumptions on $(\bx_t,Y_t)$ given in their Theorems 1--3.

The methodology developed in Section \ref{sec3} and its extension
outlined in
Section~\ref{newsec61} suggest an alternative approach to comparing
econometric forecasts. As in~(\ref{n64}), we consider the average
score difference
%
%
\begin{equation}\label{n67}
\hdel_n=n^{-1}\sum_{i=1}^n\{L(Y_t,\hY_t')-L(Y_t,\hY_t'')\},
\end{equation}
in which $\hY_t'$ and $\hY_t''$ are forecasts that are
$\g_{t-1}$-measurable. Since $k$-step ahead forecasts of $Y_t$ are
$\g_{t-1}$-measurable for any $k\geq1$, the theory applies to all
$k$-step ahead forecasts of $Y_t$, as illustrated in Table \ref{tableBS}.
Instead of hypothesis testing, our approach is targeted toward
estimating
%
%
\begin{equation}\label{n68}
\Delta_n=n^{-1}\sum_{i=1}^nE\{L(Y_t,\hY_t')-L(Y_t,\hY_t'')|\g
_{t-1}\},
\end{equation}
in which $E$ is with respect to the actual but unknown probability
measure $P$. Analogously to Theorem \ref{thm2}(i) for the case of binary
$Y_i$, we can\vspace*{2pt} apply the martingale central limit theorem to establish
the asymptotic normality of $\hdel_n-\Delta_n$. In many applications,
one can make use of the bucket structure of the type in Section \ref
{sec33} to
estimate the asymptotic variance of $\hdel_n$. In particular, this
structure is inherent in dynamic panel data in econometrics and
longitudinal data in epidemiology, which is beyond the scope of this
paper on forecasting probabilities of events and will be treated
elsewhere. Note that the bucket structure is only used in estimating
the asymptotic variance of $\hdel_n$ by (\ref{n62}), and that
Theorem \ref{thm5} and its extension outlined in Section \ref{newsec61}
imply that the
variance estimate tends to be conservative if the assumed bucket
structure actually fails to hold.

\section{Discussion}\label{sec6}

The average score $n^{-1}\sum_{i=1}^nL(Y_i,\hp_i)$ measures the
divergence of the predicted probabilities $\hp_i$, which lie between 0
and 1, from the indicator variables $Y_i$ that can only have values 0
or 1. As noted by \citet{licht}, this tends to encourage more
aggressive bets on the binary outcomes, rather than the forecaster's
estimates of the event probabilities. For example, an estimate of 95\%
probability may lead to a probability forecast of 100\% for a higher
reward associated with the indicator variable $Y_i$; see also
\citet{mason}, who gives an example in which a forecaster is encouraged
to give such ``dishonest'' forecasts. This difficulty would disappear
if one uses $L$ to compare $\hp_i$ with the actual $p_i$, rather than
with the Bernoulli($p_i$) random variable $Y_i$. Because the $p_i$ are
unknown, this is not feasible and the importance of using a proper
score $L(Y_i,\hp_i)$ to evaluate a probability forecast has been
emphasized to address the issue of dishonest forecasts. In Section \ref{sec32}
we have shown that it is possible to use $L(p_i,\hp_i)-L(p_i,\hp_i')$
for comparing two forecasters and to construct confidence intervals of
the average loss difference. A key idea underlying this development is
the linear equivalent of a loss function introduced in Section \ref{sec2}.
\citet{schervish}, Section 3, has used a framework of two-decision
problems involving these loss functions to develop a method for
comparing forecasters. Our approach that considers
$L(p_i,\hp_i')-L(p_i,\hp_i'')$ can be regarded as a further step in
this direction.

As noted in Section \ref{sec35}, an important assumption underlying
statistical inference in the verification of probability forecasts in
meteorology is that the forecast-observation pairs are independent
realizations from the joint distribution of forecasts and
observations. Although Mason [(\citeyear{mason}), page 32] has pointed
out that this
assumption cannot hold ``if the verification score is calculated using
forecasts for different locations, or if both the forecasts and
observations are not independent temporally,'' not much has been done
to address this problem other than using moving-blocks bootstrap
[\citet{mason}, \citet{wilks}] because traditional
statistical inference does not
seem to provide much help in tackling more general
forecast-observation pairs. The new approach in Section \ref{sec3} can
be used
to resolve this difficulty. It uses martingale theory to allow the
forecast-observation pairs to be generated by general stochastic
systems, without the need to model the underlying system in carrying
out the inference. The treatment of spatial dependence is also covered
in Section \ref{sec35}, in which dependence of the events at $K_t$ locations
at time $t$ is encapsulated in the highly complex joint distribution
of their generating probabilities $p_{t,1},\ldots,p_{t,K_t}$, which our
approach does not need to model in performing inference on forecast
validation. Our viewpoint in forecast evaluation is that one should
try not to make unnecessary or arbitrary assumptions on the underlying
data-generating mechanism, especially in regulatory settings such as
regulatory supervision of a bank's internal ratings models of loan
default probabilities; see Section \ref{sec33} and \citet
{2008basel}. A
convenient but incorrect data-generating model that is assumed can
unduly bias the comparison.

\section*{Acknowledgment}

The author thanks the Epidemiology unit at INSERM, France,
and the Statistics Department at Tel Aviv University for their generous
hospitality while working on revisions of this paper.


%

\printaddresses

\end{document}